\documentclass[12pt]{article}

\usepackage{amsmath}
\usepackage{amssymb}
\usepackage{latexsym}

\usepackage{enumerate}

\newtheorem{theorem}{Theorem}[section]
\newtheorem{proposition}[theorem]{Proposition}

\newtheorem{lemma}[theorem]{Lemma}

\newenvironment{proof}{\noindent{\sc Proof.}}{\quad\qed\medskip}

\newcommand{\Z}{{\mathbb Z}}

\newcommand{\Hom}{{\rm Hom}}

\newcommand{\dospunts}{\hspace{-0.1cm}:}
\newcommand{\qed}{\quad\lower0.05cm\hbox{$\Box$}}

\newcommand{\E}{{\cal E}}
\newcommand{\D}{{\cal D}}
\newcommand{\Fr}{{\rm Fr}}
\newcommand{\F}{{\cal F}}

\newcommand{\downarrowright}[1]{\downarrow
\rlap{\raise0.1cm\hbox{$\scriptstyle{#1}$}}}

\newcommand{\downarrowleft}[1]{\rlap{\kern-0.2cm
\raise0.1cm\hbox{$\scriptstyle{#1}$}}\downarrow}

\newcommand{\uparrowright}[1]{\uparrow
\rlap{\lower0.1cm\hbox{$\scriptstyle{#1}$}}}

\newcommand{\uparrowleft}[1]{\rlap{\kern-0.2cm
\lower0.1cm\hbox{$\scriptstyle{#1}$}}\uparrow}

\newcommand{\da}{\downarrow}

\newcommand{\flecha}[1]{\stackrel{#1}{\longrightarrow}}

\newcommand{\map}{{\rm map}}

\newcommand{\clas}[1]{\langle{#1}\rangle}

\newcommand{\padics}{\Z^{\raise0.05cm\hbox{\scriptsize $\wedge$}}_p}
\newcommand{\p}{^{\raise0.05cm\hbox{\scriptsize $\wedge$}}_p}
\newcommand{\3}{^{\raise0.05cm\hbox{\scriptsize $\wedge$}}_3}

\begin{document}
\title{\vspace*{0cm} Homology equivalences inducing an
epimorphism on the fundamental group and Quillen's plus construction}

\author{
{\sc Jos\'e L. Rodr\'{\i}guez and Dirk Scevenels}
\thanks{The first author was supported by DGES grant PB97-0202.\newline
}}

\date{\empty}
\maketitle

\begin{abstract}
\noindent
Quillen's plus construction is a topological construction
that kills the maximal perfect subgroup of the fundamental
group of a space without changing the integral homology
of the space. In this paper we show that there is a
topological construction that, while leaving the integral
homology of a space unaltered, kills even the intersection
of the transfinite lower central series of its fundamental
group.
Moreover, we show that this is the maximal subgroup that
can be factored out of the fundamental group without
changing the integral homology of a space.
\end{abstract}

\setcounter{section}{-1}

\section{Introduction}
As explained in \cite{Dro96}, \cite{DOS89}, Bousfield's
$H\Z$-localization $X_{H\Z}$ of a space $X$ (\cite{Bou77a})
is homotopy equivalent to its localization with respect to
a map $Bf \dospunts B \F_1 \to B \F_2$ induced by a certain
homomorphism $f \dospunts \F_1 \to \F_2$ between
free groups. This means that a space $X$ is $H\Z$-local if and
only if the induced map
$Bf^* \dospunts \map(B \F_2, X) \to \map(B \F_1, X)$
is a weak homotopy equivalence. Moreover, the effect of
$Bf$-localization on the fundamental group produces precisely
the group-theoretical $H\Z$-localization (i.e. $f$-localization)
of the fundamental group, i.e.
$\pi_1 L_{Bf} X \cong L_f(\pi_1X) \cong (\pi_1X)_{H\Z}$ for all
spaces~$X$.

%\pagebreak

A universal acyclic space for $H\Z$-localization
(i.e. $Bf$-localization), in the sense of Bousfield (\cite{Bou97}),
was studied by Berrick and Casacuberta in \cite{BC98}. They show
that nullification with respect to such a universal acyclic space
coincides with Quillen's plus construction $X^+$ for a space $X$.
Moreover, this universal acyclic space can be taken to be a
two-dimensional Eilenberg--Mac~Lane space $K(A(f),1)$, where $A(f)$
is a locally free, universal $f$-acyclic group, in the sense of
\cite{RS98}. The effect of $K(A(f),1)$-nullification on the
fundamental group is precisely given by its $A(f)$-nullification.
This $A(f)$-nullification factors out the perfect radical of a
group, i.e. the maximal perfect subgroup, which can be
obtained as the intersection of the (transfinite) derived series
of the group (\cite{CRS98}).

The intersection $\Gamma G$ of the transfinite lower central series
of a group~$G$ is also a radical, and, as observed by Bousfield
(\cite{Bou77b}), it is in fact the maximal $G$-perfect normal subgroup
of~$G$ (where a normal subgroup $H$ of $G$ is called
$G$-perfect if  $H = [H,G]$).
As explained in~\cite{CRS98}, there is an epireflection
that corresponds to this radical $\Gamma G$. By adapting the methods
used by Berrick and Casacuberta in \cite{BC98}, we describe
in Theorem~2.5 an epimorphism $g$ such that $L_gG \cong G/\Gamma G$
for all groups $G$.

We further construct a localization functor of topological spaces
that is intermediate between Quillen's plus construction and
Bousfield's homological localization, and induces localization with
respect to $g$ on the fundamental group. More precisely, we describe
in Theorem~2.9 a map $\varphi$ such that $\varphi$-localization of spaces
factors out the intersection of the transfinite lower central series of the
fundamental group, while preserving the integral homology of the space.

We furthermore show that this is in some sense the best possible
result. Indeed, we show in Proposition~3.2 that the maximal subgroup
that can be factored out of the fundamental group without altering the
integral homology of a space, is precisely the intersection of the
transfinite lower central series of the fundamental group.

Finally, we turn our attention to the question whether there exists a
localization functor of spaces which induces localization with respect
to a universal epimorphic $H\Z$-equivalence of groups (in the sense
of \cite{RS98}) on the fundamental group, without changing the integral
homology of a space. In fact, we show in Proposition~3.3 that the answer
to this question is affirmative if and only if, for all groups $G$, the
kernel of the $H\Z$-localization homomorphism $G \to G_{H\Z}$ coincides
with $\Gamma G$.

%%\bigskip\noindent
%%{\bf Acknowledgements.}
%%We are indebted to Carles Casacuberta for his suggestions and helpful
%%comments.
%%Discussions with Ran Levi and Carles Broto are warmly appreciated.
%%The second named author also wishes to thank the Centre de
%%Recerca Matem\`atica for its hospitality during his stay in which this
%%work was completed.

\setcounter{equation}{0}
\section{Preliminaries}
For the convenience of the reader, we recall some terminology and bacis
facts about localization with respect to a given continuous map
(see e.g. \cite{Bou97}, \cite{Dro96}).
Given a map $f \dospunts A \to B$,
a space $X$ is called {\it $f$-local\/} if the induced map
\[
f^* \dospunts \map(B,X) \to \map(A,X)
\]
is a weak homotopy equivalence.
For every space $X$ there is a map $l_X \dospunts X \to L_fX$,
which is initial among all maps from $X$ into $f$-local spaces.
$L_f$ is called the localization functor with respect to~$f$.
A map $\phi$ is called an
{\it $f$-equivalence\/} if $L_f\phi$ is a
homotopy equivalence.
Further a space $X$ is called {\it $f$-acyclic\/} if $L_fX$ is
contractible.
In the special case where $f$ is of the form $f \dospunts A \to *$,
the $f$-localization of a space~$X$ is also denoted by
$P_AX$ and it is called the {\it $A$-nullification\/} of $X$.
The {\it localization class\/} of $f$, denoted by $\clas{f}$, is
defined as the collection of all maps~$g$ such that $L_g $ is
naturally equivalent to~$L_f$.
When $f$ is of the form $f \dospunts A \to *$, the class of $f$
is simply denoted by $\clas{A}$ and is called the
{\it nullification class\/} of $A$.
One says that $\clas{f} \leq \clas{g}$ if and only if
there is a natural transformation of localization functors
$L_f \to L_g$. Recall that this is equivalent to every $g$-local
space being $f$-local, to every $f$-equivalence being a
$g$-equivalence, or to $f$ being a $g$-equivalence.
As shown by Bousfield in \cite{Bou97}, the collection of
localization functors with respect to maps is a small-complete
lattice for this partial order relation.
Note that a space $X$ is $f$-acyclic if and only if
$\clas{X} \leq \clas{f}$.
Furthermore, in \cite{Bou97} Bousfield proved that every
localization class has a best possible approximation by a
nullification class. More precisely, for any map~$f$, there is a
maximal nullification class $\clas{A(f)}$ which is smaller than
$\clas{f}$ in the lattice of localization classes. The space $A(f)$
is called a {\it universal $f$-acyclic space}, since a space $X$ is
$f$-acyclic if and only if $X$ is $A(f)$-acyclic.

An important algebraic tool in studying localization with respect
to a given map is given by its discrete analogue in the category
of groups (see e.g. \cite{Bou97}, \cite{Cas95}).
For a given group homomorphism $f \dospunts A \to B$, a group $G$
is called {\it $f$-local\/} if the induced map of sets
\[
f^* \dospunts \Hom(B,G) \to \Hom(A,G)
\]
is a bijection. For every group $G$ there is a homomorphism $l_G
\dospunts G \to L_fG$, which is initial among all homomorphisms from
$G$ into $f$-local groups. This allows to introduce the localization
functor $L_f$ with respect to $f$. In an obvious way, one speaks of
$f$-equivalences and $f$-acyclic groups, and, if $f$ is of the form
$f \dospunts A \to 1$,  of $A$-nullification. It makes also sense to
speak of the localization class of a homomorphism $f$, denoted by
$\clas{f}$, which in the special case when $f$ is of the form $f
\dospunts A \to 1$, is simply denoted by $\clas{A}$ and is called the
nullification class of $A$. Furthermore, the collection of all
localization classes of homomorphisms is again a small-complete
lattice for an obvious partial order relation. In \cite{RS98} it was
proved that every localization class of a homomorphism has a best
possible approximation by a nullification class and by the
localization class of an epimorphism. More precisely, for any
homomorphism $f$, there is a maximal nullification class
$\clas{A(f)}$ and a maximal class $\clas{\E(f)}$ where $\E(f)$ is an
epimorphism such that
\[
\clas{A(f)} \leq \clas{\E(f)} \leq \clas{f}.
\]
The group $A(f)$ is called a {\it universal $f$-acyclic group}, since
a group $G$ is $f$-acyclic if and only if $G$ is $A(f)$-acyclic.
Accordingly, $\E(f)$ is called a {\it universal epimorphic
$f$-equivalence},  since an epimorphism $g$ is an $f$-equivalence if
and only if $g$  is an $\E(f)$-equivalence.

It was shown in \cite{RS98} (see also \cite{CRS98}) that
to any localization class $\clas{f}$ of an epimorphism~$f$, there is
associated a radical $R_f$ on the category of groups such that
$L_f G \cong G/R_fG$ for all groups $G$. (Recall that a radical
$R$ is a functor which assigns to every group $G$ a normal subgroup
$RG$ in such a way that every homomorphism $G \to K$ restricts to
$RG \to RK$ and such that $R(G/RG) = 1$.)
In fact, there is a bijective correspondence between epireflections
(i.e. idempotent functors $L$ on the category of groups for which
$G \to LG$ is an epimorphism for any group $G$) and radicals.
Furthermore, if the class $\clas{f}$ is actually a nullification
class, then the associated radical $R_f$ is idempotent, meaning that
$R_fR_fG = R_fG$ for all groups $G$.

\setcounter{equation}{0}
\section{A universal epimorphic $H \Z$-map}
It is well known that Bousfield's $H\Z$-localization $X_{H\Z}$ of a
space $X$ (\cite{Bou77a}) is homotopy equivalent to its localization
with respect to a map $Bf \dospunts B \F_1 \to B \F_2$ induced by a
certain homomorphism $f \dospunts \F_1 \to \F_2$ between free groups
(\cite{Dro96}, \cite{DOS89}). In fact, the homomorphism~$f$ can be
taken to be the free product of a set of representatives of
isomorphism classes of homomorphisms between countable, free groups
inducing an isomorphism on the first integral homology group.
Furthermore, the effect of $Bf$-localization on the fundamental group
is to produce its $f$-localization, which coincides with the
group-theoretical $H\Z$-localization of the fundamental group. In
\cite[Proposition~4.2]{BC98} the authors show that a universal
$f$-acyclic group~$A(f)$ (in the sense of \cite{RS98}) can be taken
to be the free product of a set of representatives of all isomorphism
classes of countable, locally free, perfect groups. The key lemma
here is a result due to Heller (\cite{Hel80}), which states that for
every element~$x$ in any perfect group $P$, there exists a countable,
locally free, perfect group~$D$ and a homomorphism $D \to P$
containing $x$ in its image.
We can similarly show the following lemma.
%We first prove a \lq \lq relative"
%version of \cite[Lemma 5.7]{Hel80}.
Recall that a normal subgroup $H$
of a group $G$ is called {\it $G$-perfect\/} if  $H = [H,G]$.
In fact, any group $G$ has a maximal $G$-perfect
subgroup, which we denote by~$\Gamma G$, and that can be obtained as the
intersection of the transfinite lower central series of~$G$ (see \cite{Bou77b}).

\begin{lemma}
Let $A$ be any group with an $A$-perfect normal subgroup $K$.
Let $z$ be any element in $K$. Then there exist a countable,
locally free group $D$ and a homomorphism $\psi \dospunts D \to A$ such that $z$ belongs
to the image of some element in $\Gamma D$.
\qed
\end{lemma}
\begin{proof}
Since $K=[K,A]$ we can express $z$ as a finite product of commutators
$z=\Pi_i[k_i,a_i]$ where $k_i\in K$, $a_i\in A$.
Now again, for each $k_i$, we can find finitely many
elements $k_{i,l}\in K$, $a_{i,l}\in A$,
such that $k_i=\Pi_l[k_{i,l},a_{i,l}]$.
Similarly for $k_{i,l}$ we have $k_{i,l}=\Pi_m[k_{i,l,m},a_{i,l,m}]$, with
$k_{i,l,m}\in K$ and $a_{i,l,m}\in A$, etc.
 For each element appearing in these commutator decompositions
we choose a free generator, which we denote by $x$, $x_i$, $y_i$, $x_{i,l}$
$y_{i,l}$, $x_{i,l,m}$, $y_{i,l,m}$, etc. Then we define $D$ as the following
direct limit of free groups:
$$
\Fr(x) \stackrel{\psi_0}{\longrightarrow}
\Fr(x_{i}, y_{i}) \stackrel{\psi_1}{\longrightarrow}
\Fr(x_{i,l}, y_{i,l}, y_i)  \stackrel{\psi_2}{\longrightarrow}
\Fr(x_{i,l,m}, y_{i,l,m}, y_{i,l}, y_i)
\stackrel{\psi_3}{\longrightarrow}\cdots
$$
where $\psi_j$ are defined as the decompositions above, that is
$\psi_0(x)= \Pi_i[x_i,y_i]$; $\psi_1(x_i)=\Pi_l[x_{i,l},y_{i,l}]$,
and the identity on the $y_i$'s;
$\psi_2(x_{i,l})= \Pi_m[x_{i,l,m},y_{i,l,m}]$,
and the identity on the other generators, etc.
Observe that $D$ is countable and locally free.
The obvious homomorphism $\psi:D\to A$ is defined as
$\psi(x)=z$, $\psi(x_i)=k_i$, $\psi(y_i)=a_i$,
$\psi(x_{i,l})=k_{i,l}$, $\psi(y_{i,l})=a_{i,l}$, etc.
It is clear by construction that $x\in \Gamma D$ and $\psi(x)=z$ as required.
\end{proof}

%Lemma 2.1 enables us to prove the following result.

\begin{proposition}
Let $G$ be any group. Then the following assertions are equivalent:
\begin{itemize}
\item[\rm (i)] For every group $A$ and every $A$-perfect normal
subgroup $K$ of~$A$, the restriction of any homomorphism
$A \to G$ to $K$ is trivial;
\item[\rm (ii)] For every countable, locally free group $A$ and every
$A$-perfect normal subgroup $K$ of~$A$, the restriction of any
homomorphism $A \to G$ to $K$ is trivial;
\item[\rm (iii)] For every countable, locally free group $A$,
the restriction of any homomorphism $A \to G$ to $\Gamma A$ is trivial.
\qed
\end{itemize}
\end{proposition}

The subgroup $\Gamma G$ actually defines a radical on the category of
groups, and, hence, by~\cite{CRS98}, the assignation $G \to LG = G/\Gamma G$
is an epireflection. Our aim is to show that this epireflection is
singly generated. More precisely, we want to exhibit a homomorphism $g$ such
that $L_gG \cong LG = G/\Gamma G$ for any group $G$. Moreover, by \cite{RS98}
we know that it is possible to choose $g$ to be an epimorphism.

Let $\Lambda$ be a set of representatives of
isomorphism classes of countable, locally free
groups. If $\D'$ denotes the
normal closure of $\bigcup_{\Lambda} \Gamma D$ in the free product
$\D$ of all groups~$D$ in~$\Lambda$, which we denote by
$\Fr_{\Lambda} D$, then we have a short exact sequence
(cf.~\cite[Exercise~6.2.5]{Rob96})
\begin{equation}
\label{ses1}
\D' \rightarrowtail \D = \Fr_{\Lambda} D \twoheadrightarrow
\D/\D' \cong \Fr_{\Lambda}(D/\Gamma D),
\end{equation}
where $\Fr_{\Lambda}(D/\Gamma D)$ denotes the free product of all
groups $D/\Gamma D$ for which (the isomorphism class of)
$D$ belongs to $\Lambda$.

\begin{proposition}
Localization with respect to the epimorphism
\[
h \dospunts \D \twoheadrightarrow \D/\D'
\]
given in (\ref{ses1}) satisfies $L_h G \cong G/ \Gamma G$ for
all groups $G$.
\end{proposition}
\begin{proof}
Observe that a group $G$ is $f$-local for any given epimorphism
$f \dospunts A \twoheadrightarrow B$ if and only if the restriction
of any homomorphism $A \to G$ to the kernel of $f$ is trivial.
The proof is now completed by using Proposition~2.2.
\end{proof}

%\pagebreak

Note that we can partition $\Lambda$ into $\Lambda_1$, containing
all the representatives of the isomorphism classes of countable,
locally free, perfect groups, and its complement $\Lambda_1^c$.
We then can write
\begin{eqnarray*}
\clas{h} & = & \clas{\Fr_{\Lambda_1} (D \to D/\Gamma D)}
     * \clas{\Fr_{\Lambda_1^c} (D \to D/\Gamma D)} \\
         & = & \clas{\Fr_{\Lambda_1}(D \to 1)} *
    \clas{\Fr_{\Lambda_1^c} (D \to D/\Gamma D)} \\
         & = & \clas{\F} *
               \clas{\Fr_{\Lambda_1^c} (D \to D/\Gamma D)},
\end{eqnarray*}
where $\F$ is the universal $H\Z$-acyclic group defined in~\cite{BC98}.
(Here we denote by $\clas{f_1} * \clas{f_2}$ the least upper
bound of the classes $\clas{f_1}$ and $\clas{f_2}$ in the lattice of
localization classes, and we have used the fact that the free product
$f_1*f_2$ is a representative of this least upper bound.)

However, it is possible to give another description of $\clas{h}$,
which will be more useful later on.
Indeed, note that,
if $f_1 \dospunts A \twoheadrightarrow B$ is an epimorphism,
and $f_2 \dospunts B \to C$ is any homomorphism, then
$\clas{f_2 \circ f_1} = \clas{f_1 * f_2} = \clas{f_1} * \clas{f_2}$.
This enables us to prove the following preliminary result.

\begin{lemma}
Let $f \dospunts A \twoheadrightarrow B$ be any epimorphism.
Then $\clas{f} = \clas{l_A}$,
where $l_A \dospunts A \to L_fA$ denotes the
localization homomorphism.
\end{lemma}
%\begin{proof}
\noindent
{\sc Proof.}
Since $f$ is an epimorphism, we infer from \cite[Theorem~2.1]{RS98} that
$l_A$ is an epimorphism.
Hence,
\[
\clas{f} = \clas{f} * \clas{l_B} = \clas{l_B \circ f} =
\clas{L_ff \circ l_A} = \clas{L_ff} * \clas{l_A} = \clas{l_A}.\qed
\]
%\end{proof}

The above lemma now proves the following alternative description
of $\clas{h}$.

\begin{theorem}
Consider the natural homomorphism $g \dospunts \D \to \D/ \Gamma \D$,
where $\D$ is as defined in (\ref{ses1}). Then
$\clas{g} = \clas{h}$, where $h$ is defined in (\ref{ses1}).
In other words, $L_g G \cong G/ \Gamma G$ for all
groups $G$.
\qed
\end{theorem}

In fact, the homomorphism $g$ is a \lq \lq universal epimorphic
$H\Z$-map", as we next show.
Recall from \cite{Bou77b} that a group homomorphism~$g$ is called an
{\it $H\Z$-map\/} if $H_1g$ (i.e. the homomorphism induced by $g$
on the first integral homology group) is an isomorphism and $H_2g$
is an epimorphism.
We first need  a characterization of the epimorphisms that are
$H\Z$-maps (cf.~\cite{Bou77b}).

\begin{lemma}
Let $h$ be an epimorphism $h \dospunts A \twoheadrightarrow B$
with kernel $K$. Then $h$ is an $H\Z$-map if and only if $K$ is
$A$-perfect.
\end{lemma}
\noindent
{\sc Proof.}
This is an obvious consequence of the 5-term exact sequence
\[
H_2(A) \to H_2(B) \to K/[K,A] \to H_1(A) \to H_1(B) \to 0. \qed
\]

\begin{proposition}
Let $G$ be any group and let $g \dospunts \D \to \D/\Gamma \D$,
where $\D$ is as defined in (\ref{ses1}).
The homomorphism $ G \to L_gG \cong G/\Gamma G$ is terminal among all
epimorphic $H\Z$-maps going out of $G$.
\qed
\end{proposition}

Observe that the epimorphism $g$ that we have constructed is not a
universal epimorphic $H\Z$-equivalence (in the sense of \cite{RS98}).
Indeed, there are \lq \lq more"
epimorphic $H\Z$-equivalences than epimorphic $H\Z$-maps.
However, if we denote by $A(f)$, resp. $\E(f)$ a universal $H\Z$-acyclic
group, resp. a universal epimorphic $H\Z$-equivalence, then there are natural
homomorphisms
\[
G \to P_{A(f)}G \to L_g G \cong G/\Gamma G \to L_{\E(f)}G \to G_{H\Z},
\]
for any group $G$,
where $G \to L_{\E(f)}G \to G_{H\Z}$ is an epi-mono factorization.
Moreover, for many groups $G$ (e.g. finite groups, nilpotent groups, or
more generally, groups for which the lower central series stabilizes),
we have isomorphisms
$L_gG \cong G/\Gamma G \cong L_{\E(f)}G \cong G_{H\Z}$
(cf. \cite{Bou77b}).

We now want to realize localization with respect to
$g \dospunts \D \twoheadrightarrow \D/\Gamma \D$ topologically,
by exhibiting a localization functor of spaces that induces
$g$-localization on the fundamental group and which is intermediate
between Quillen's plus construction and Bousfield's $H\Z$-localization
(and, hence, does not change the integral homology of a space).

We will need the following proposition, which is similar to results
obtained in~\cite{Cas95} and~\cite{CRT98}.

\begin{proposition}
Let $\psi \dospunts A \to B$ be any map which induces an epimorphism
$\psi_*  = \pi_1(\psi) \dospunts \pi_1(A) \twoheadrightarrow \pi_1(B)$
and suppose that $A$ is a CW-complex of dimension at most two.
Then $\psi$-localization of spaces is $\pi_1$-compatible, i.e.
\[
\pi_1(L_\psi X) \cong L_{\psi_*}(\pi_1X),
\]
for all spaces $X$.
\end{proposition}
\begin{proof}
For any space $X$, the map $X \to L_\psi X$ is a $\psi$-equivalence.
Hence, by \cite[Proposition~3.3]{Cas95}, the induced homomorphism
$\pi_1(X) \to \pi_1(L_\psi X)$ is a $\psi_*$-equivalence.
Moreover, we claim that $\pi_1(L_\psi X)$ is $\psi_*$-local.
To see this, it suffices to prove that the restriction of any
homomorphism $\ell \dospunts \pi_1A \to \pi_1(L_\psi X)$ to
$\ker \psi_*$ is trivial. However, since the dimension of $A$ is
at most two,
there exists a map $\xi \dospunts A \to L_\psi X$ inducing $\ell$
on the fundamental group. Since $L_\psi X$ is $\psi$-local, we infer
that there exists a map $\chi \dospunts B \to L_\psi X$ such that
$\chi \circ \psi \simeq \xi$, which implies that the restriction
$\ell | \ker \psi_* = (\pi_1(\chi) \circ \psi_*) | \ker \psi_* $ is
trivial.
\end{proof}

Now choose a two-dimensional CW-complex $M\D$ such that
$\pi_1 M\D = \D$. Attach 2-cells to $M\D$, thereby
obtaining a map $i \dospunts M\D \hookrightarrow C$ which
induces
\[
g \dospunts \D \twoheadrightarrow \D/\Gamma \D
\]
on the
fundamental group. We then obtain a diagram
\begin{eqnarray}
\begin{array}{ccccccccc}
  &  & \pi_2 M\D & \to & \pi_2 C & &  & & \\
  &  & \da       &     & \da     & &  & & \\
0 & \to & H_2 M\D & \to & H_2 C & \to & H_2(C,M\D) & \to & 0 \\
  &     &  \da    &     & \da   &     &            &     &   \\
  &     & H_2(\pi_1 M\D) & \to & H_2(\pi_1 C) & \to & 0 & & \\
  &     &  \da           &     &  \da         &     &   & & \\
  &     &  0             &     &  0           &     &   & & \\
\end{array}
\nonumber
\end{eqnarray}
Moreover, since $H_2(\pi_1 C) = H_2(\D/\Gamma \D) = 0$, we infer
that $H_2 C$ is an epimorphic image of $\pi_2 C$. This means that
we can kill $H_2 C$ by attaching 3-cells to $C$, through a map
$j \dospunts C \to C'$. It is now easily verified that the
composition $\varphi$ of
\begin{equation}
\label{varphi}
M\D \flecha{i} C \flecha{j} C'
\end{equation}
is an integral homology equivalence and that $\varphi$ induces the
homomorphism
$g \dospunts \D \twoheadrightarrow \D/\Gamma \D$ on the
fundamental group (cf. \cite[Lemma~6.1]{Bou77a}).

\begin{theorem}
Let $\varphi \dospunts M\D \to C'$ be the composition
given in (\ref{varphi}).
Then $\varphi$-localization of spaces
induces $g$-localization on the fundamental group.
Moreover, for any space $X$, there are natural maps
\[
 X \to X^+ \to L_{\varphi} X \to X_{H\Z}.
\]
\end{theorem}

%\pagebreak

\begin{proof}
Since $M\D$ is a two-dimensional CW-complex and since $\varphi$
induces an epimorphism on the fundamental group, we know by
Proposition~2.8 that $\varphi$-localization is $\pi_1$-compatible,
so that it induces $g$-localization on the fundamental group.
The second claim is obvious, since $\varphi$ is clearly an
integral homology equivalence, and the perfect radical of $\pi_1(L_\varphi X)$
being trivial (cf.~\cite{BC98}).
\end{proof}

To see that the natural maps given in Theorem~2.9 are not equivalences
in general, and thus that we have constructed a functor which is really
different from Quillen's plus construction and from $H\Z$-localization,
observe the following.
On one hand,
\[
L_{\varphi}(S^1 \vee S^1) \simeq S^1 \vee S^1 \simeq (S^1 \vee S^1)^+.
\]
Indeed, the fact that $ \Z * \Z$ is $g$-local implies that
$S^1 \vee S^1$ is $\varphi$-local.
On the other hand,
$\pi_1(L_\varphi B\Sigma_3) \cong L_g(\pi_1 B\Sigma_3) \cong \Z/2$,
while
$\pi_1 B\Sigma_3^+ \cong P_{A(f)}(\pi_1 B\Sigma_3) \cong \Sigma_3$.

\setcounter{equation}{0}
\section{Homology equivalences inducing an epimorphism on
the fundamental group}
In this section we want to explore some immediate consequences of our
results. In particular, we want to show that there are some restrictions
on (integral) homology equivalences that induce an epimorphism on the
fundamental group.

\begin{proposition}
Let $\psi \dospunts X \to Y$ be an integral homology equivalence of spaces
such that the induced homomorphism $f = \pi_1 \psi$ is an epimorphism.
Then $R_fG \subset \Gamma G$ for all groups~$G$, where $R_f$ denotes
the radical associated to the epireflection class $\clas{f}$.
\end{proposition}
\begin{proof}
By hypothesis we know that $f$ is an epimorphic $H\Z$-map. Hence,
$\clas{f} \leq \clas{g}$, where $g$ denotes the universal epimorphic
$H\Z$-map of Theorem~2.5. Hence, there  are natural homomorphisms $G
\twoheadrightarrow L_fG \cong G/R_fG \twoheadrightarrow L_gG \cong
G/\Gamma G$, for any  group~$G$.
\end{proof}

In particular, we have the following result.

\begin{proposition}
Let $\psi \dospunts X \to Y$ be an integral homology equivalence of spaces
such
that the induced homomorphism $f = \pi_1 \psi$ is an epimorphism.
Then $\ker f \subset \Gamma \pi_1 X$.
\end{proposition}
\begin{proof}
Since $f$ is an epimorphic $H\Z$-map, we know that $\ker f$ is
$\pi_1 X$-perfect.
\end{proof}

%\pagebreak

In other words, for any space $X$, the maximal subgroup
that can be factored out of $\pi_1X$ without altering the integral
homology of $X$, is precisely $\Gamma \pi_1 X$.
In particular, this implies that there is a restriction on the
possibility of realizing topologically a universal epimorphic
$H\Z$-equivalence of groups (i.e. of finding an
integral homology equivalence of spaces which induces localization with
respect
to a universal epimorphic $H\Z$-equivalence of groups
on the fundamental group).

\begin{proposition}
Let $\E(f)$ denote a universal epimorphic $H\Z$-equivalence of groups.
Then the following assertions are equivalent:
\begin{itemize}
\item[(i)] There exists an integral homology equivalence $\psi \dospunts X
\to Y$ of spaces such that $\pi_1 \psi = \E(f)$ and $L_\psi$ is
$\pi_1$-compatible;
\item[(ii)] $\ker(l_{H\Z} \dospunts G \to G_{H\Z}) = \Gamma G$
for all groups $G$.
\end{itemize}
\end{proposition}
\begin{proof}
To see that (i) implies (ii), it suffices to show that
$\ker l_{H\Z} \subset \Gamma G$.
However, this is an immediate consequence of Proposition~3.2, since
for every group $G$ we have
\[
L_{\pi_1 \psi} G \cong L_{\E(f)} G \cong G/\ker l_{H\Z}.
\]
Finally (ii) implies (i), as is shown by our construction of
$\varphi$ in (\ref{varphi}).
\end{proof}

\setlength{\baselineskip}{0.5cm}

\vskip 0.5 cm

\setlength{\baselineskip}{0.6cm}

\bigskip\noindent
\'{A}rea de Geometr\'{\i}a y Topolog\'{\i}a, CITE III, Universidad de Almer\'{\i}a\newline
E--04120 La Ca\~{n}ada de San Urbano, Almer\'{\i}a, Spain \newline e-mail: {\tt jlrodri@ual.es}

\bigskip\noindent
Departement Wiskunde, Katholieke Universiteit Leuven \newline
Celestijnenlaan 200 B, B--3001 Heverlee, Belgium \newline
e-mail: {\tt dirk.scevenels@wis.kuleuven.ac.be}

\end{document}